\documentstyle[proceedings,amssymb,numreferences]{crckapb}

\begin{opening}
\title{THE FUKAYA TYPE CATEGORIES\protect\\
       FOR ASSOCIATIVE ALGEBRAS}
\subtitle{}

\author{R.NEST}
\institute{University of Copenhagen\\
           2100 Copenhagen, Denmark}
\author{B.TSYGAN}
\institute{Penn State University\\
           University Park, PA, 16803, USA}

\end{opening}

\runningtitle{THE FUKAYA TYPE CATEGORIES}

\begin{document}

% The \begin{document} command comes after the \end{opening}
% command.

\def\pf{\hfill $\blacksquare$}
\def\ub{\underbar}
\def\c{\cite}
\def\fr{\frac}
\def\fg{\goth{g}}
\def\Z{{\Bbb Z}}
\def\W{\Bbb W}
\def\C{{\Bbb C}}
\def\call{\cal L}
\def\ll{\Bbb L}
\def\n{\nabla}
\def\uend{\underline{End} }
\def\alt{Alt}
\def\ext{Ext}
\def\pf{\Box}
\def\hom{Hom}
\def\iff{iff}
\def\lie{Lie}
\def\tor{Tor}
\def\per{per}
\def\g{\frak{g}}
\def\endo{endomorphisms }
\def\hoc{Hochschild }
\def\st{such that }
\def\asa{(A,A)}
\def\wtd{\widehat{D}}
\def\wdd{\widetilde{D}}
\def\oa{\overline{A}}
\def\ud{\underline{D}}
\def\ue{\underline{E}}
\def\ux{\underline{x}}
\def\unn{\underline{\n}}
\def\sd{D\hskip-.09in{/}\,}%\not{D}}
\def\ssd{\not{D}}
\def\usd{\underline{\ssd}}
\def\rc{\overline{C}^{\lambda}}
\def\rca{\rc_{*-1}(A)}
\def\hc{\overline{HC}}
\def\eH{\cal{H}}
\def\l{\frak{h}}
\def\glk{\frak{gl}(k)}
\def\bimod{_{\frak{gl}(k[\eta])}U(\frak{gl}(A[\epsilon,\eta]))_{\frak{gl}(k[\epsilon])}}
\def\rightisoarrow{\;\widetilde{\rightarrow}\;}
\def\H{H}
\def\Coker{{\em Coker}}
\def\Ker{{\em Ker}}
\def\Im{{\em Im}}
\def\deg{{\em deg}\;}
\def\ccp{CC_*^{per}}
\def\ah{{\Bbb A}^{\hbar}}
\def\ahc{{\Bbb A}^{\hbar}_c}
\def\E{{\cal E}_A^*}
\def\hb{\hbar}
\def\*{\star}
\def\cm{C^{\infty}(M)[[\hbar]]}
\def\ccm{C_c^{\infty}(M)[[\hbar]]}
\def\Tr{Tr}
\def\tr{tr}
\def\Cc{\check{C}}
\def\varksi{\xi}
\def\R{\Bbb R}
\def\W{\Bbb W}
\def\wu{W^*(\goth{u}(n)}
\def\wh{W^*(\goth{h}}
\def\Nabla{D}
\newtheorem{th}{Theorem}
\newtheorem{definition}[th]{Definition}
\newtheorem{proposition}[th]{Proposition}
\newtheorem{lemma}[th]{Lemma}
\newtheorem{question}[th]{Question}
\newtheorem{corollary}[th]{Corollary}
\newtheorem{remark}[th]{Remark}
\newtheorem{example}[th]{Example}
\newtheorem{assumption}[th]{Assumption}
\newtheorem{notation}[th]{Notation}
\newtheorem{conjecture}[th]{Conjecture}
%

%%%%%%%%%%%%%
\section{Introduction} \label{intro}
In his study of Lagrangian intersections Fukaya discovered a new algebraic structure on Morse-Floer complexes of loop (path) spaces of symplectic manifolds (\c{Fu}). He constructed a ``category'' whose objects are symplectomorphisms and whose space of morphisms are Morse-Floer complexes. The composition in this ``category'' yields multiplication on Floer cohomology. At the level of cochains this composition is ``strongly associative up to homotopy''; in other words, it defines a structure of an $A_{\infty}$ category.

It turns out that a very similar algebraic structure can be defined in terms of noncommutative differential geometry for any associative ring. When the ring in question is a result of deformation quantization of a symplectic manifold then the resemblance with the Fukaya construction becomes especially strong.

 It is well known that most geometric properties of a manifold can be recovered from the algebra $A$ of functions on this manifold.  By (affine) noncommutative geometry one usually means extending corresponding constructions to the case when $A$ is not necessarily commutative.  For example, 
consider the standard cochain complex $C^*(A,A)$ computing the Hochschild cohomology 
$$
        H^*(A,A)=Ext^*_{A\otimes A^o}(A,A)
$$
 and the standard chain complex computing the Hochschild homology
$$
          HH_*(A)=H_*(A,A)=Tor_*^{A\otimes A^o}(A,A).
$$
 It has been shown in \c{HKR} that, if $A$ is the algebra of regular functions on an affine nonsingular algebraic variety $X$ over a field $k$ of characteristic zero, then 
$$
     HH_i(A) \simeq \Omega^i_{X/k} 
$$
$$
      H^i(A,A)\simeq \Gamma (X, \wedge ^{i}TX)
$$
The standard cochain complex $C^*(A,A)$ is a differential graded associative algebra (with the cup product which can also be interpreted as the Yoneda product). 
On the other hand, the standard chain complex $C_*(A,A)$ is not a differential graded algebra unless $A$ is commutative.
It turns out, however,
that if one considers all differential operators on $\Omega^*(X)$
(not just  zero order multiplication operators) one gets a ring
which has a noncommutative generalization. The construction we are going
to present rests upon an idea of Yu.I.Manin (\c{M}).

Let $R$ be a monoidal category.  For an object $A$ of
$R$ one calls an object $\uend \ A$ of $R$ the inner object of  \endo
of $A$ if there are natural morphisms
\begin{equation} \label{eq:inend}
	\uend \ A\times\uend \ A\to\uend \ A\quad;\quad
	\uend \ A\times A\to A
\end{equation}
which are universal and associative in a natural way.
In the category of associative algebras inner objects of \endo
do not exist (for example, the set of \endo of an algebra is not a 
linear space).  In this paper, however, we will show that, if one takes 
for $\uend \ A$ the differential graded algebra $C^*(A,A)$ of Hochschild
cochains, the maps (\ref{eq:inend}), in a sense, still exist.  More precisely,
they exist if one passes from the category of algebras
to the category of complexes by means of some well known homological
functors.

Let $A$ be an associative unital algebra over a commutative unital
ground ring $k$.  Consider $A$ as a bimodule over itself; by $\E$ we denote the differential graded algebra
$(C^*\asa,\delta,\smile)$ which is the standard complex for computing
$\ext_{A\otimes A^{\circ}}^*\asa=H^*\asa$ (the \hoc cohomology)
equipped with the Yoneda product (Sect. \ref{hocochain}).  We construct the map of
complexes
\begin{equation}  \label{eq: oper}
	\bullet : C_*(A) \otimes C_*(\E) \to C_*(A) 
\end{equation}
\begin{equation} \label{eq: oper1}
	\bullet : C_*(\E) \otimes C_*(\E) \to C_*(\E) 
\end{equation}
Here $C_*$ stands for the Hochschild complex computing
$HH_*(A)$ or for the periodic cyclic complex.

 For the case of bar complex an analogous operation was discovered by Getzler - Jones and by Gerstenhaber - Voronov (\c{GJ}, \c{GV}).
Therefore one can say informally that in the homotopy category of
(differential graded) algebras the inner objects of \endo always exist.

Next we construct construct ``the Fukaya category'' of an associative algebra $A$; the objects of this category are automorphisms of $A$ and, for two endomorphisms $\alpha$ and $\beta$, the complex $Hom(\alpha, \beta)$ is the twisted cochain complex $C_*(\E, _{\alpha} \cal{E}^* _{\beta})$. We also construct a functor putting in correspondence to $\alpha$ the twisted chain complex $C_*(A, A_{\alpha})$. Conjecturally, these are an $A_{\infty}$ category and an $A_{\infty}$ functor.

In other words, we construct the products
\begin{equation} \label{eq:fuk}
 \bullet :C_*(A,A_{\alpha}) \otimes C_*(\E, _{\alpha}{\cal{E}}^*_{\beta}) \rightarrow  C_*(A, A_{\beta})
\end{equation}
\begin{equation} \label{eq:fuk1}
\bullet :C_*(\E,\; _{\alpha}{\cal{E}}^*_{\beta}) \otimes C_*(\E, \;_{\beta}{\cal{E}}^*_{\gamma}) \rightarrow C_*(\E, \;_{\alpha}\cal{E}^*_{\gamma})
\end{equation}
associative up to homotopy. 

Depending on whether one uses finite or infinite cochains, one can define two versions of the homology category ${\cal H}^*_A$ or  ${\cal H}^*_{A, \infty}.$

Next we consider some examples. First, in case when $A = k[V]$ where $V$ is affine nonsingular, one has 
$$ {\cal{H}}^*_A(1,1) = {\cal {D}}(\Omega^* (V))^{op}$$
where ${\cal{D}}$ is the ring of differential operators and $^{op}$ stands for the opposite ring. We consider a partial case in which we prove that 
$$ {\cal{H}}^*_A(\alpha,\alpha) = {\cal {D}}(\Omega_* (V^{\alpha}))^{op}$$
where $V^{\alpha}$ is the fixed point set; $ {\cal{H}}^*_A(\alpha,\beta)$ can be obtained from a standard $D$-module construction. On the other hand, 
$$ {\cal{H}}_{A, \infty}^*(1,1) = End(\Omega^* (V))$$
and it seems that in a rather big generality (whether in the commutative case or not)
$$ {\cal{H}}^*_{A, \infty}(\alpha,\beta) = Hom (H_*(A, A_{\alpha}), H_*(A, A_{\beta}))$$

Next we consider the case when $A$ is a deformed ring of functions on a symplectic manifold $M$ over the field of Laurent series ${\C}[\hbar, \hbar ^{-1}]]$ (\c{BBFLS}). We show that 
\begin{equation}  \label{eq:petli}
{\cal{H}}^*_A(1,1) = H^*(M^{S^1}, {\C}[\hbar, \hbar ^{-1}]]).
\end{equation}
After this we study the subcategory of the category ${\cal{H}}_A$ in a partial case when $A$ is a deformed ring of functions on $S^2$ and the automorphisms are rotations. We show how to compute $ {\cal{H}}^*_A(\alpha,\beta)$ in terms of the cohomology of the space of loops based on points $x$ for which $\alpha x = \beta x.$ 

Note the resemblance with the Fukaya construction from \c{Fu}. Recall that deformation quantization is an algebra indexed by a formal parameter from $H^2(M).$ Its construction, but not its isomorphism class, depends on an additional geometric structure, like a K\"{a}hler structure on $A.$ It looks like there is, together with the categories ${\cal{H}}^*_A$ and ${\cal{H}}^*_{A, \infty},$ some third, semi-infinite version which is more closely related to the Fukaya category and the quantum cohomology. It also seems plausible, in light of the formality conjecture of Kontsevich (\c{K}), that there exists a quantum cohomology theory for any Poisson manifold.

{\ub {Aknowledgements.}} The idea that the cohomology of deformed algebras may be related to Lagrangian intersections is due to Boris Feigin. We thank M. Kontsevich, A. Radul and A. Voronov for helpful discussions.

The authors are very grateful to the organizers of the Ascona meeting.

\section{The Hochschild cochain complex} \label{hocochain}

Let $A$ be a graded algebra with unit over a commutative unital ring
$k$. Let $M$ be a graded bimodule over $A.$  A \hoc $d$-cochain is a linear map $A^{\otimes d}\to M$.  Put,
for $d\geq 0$,
$$
	C^d (A,M) =\hom_k(\oa^{\otimes d},M)
$$
where $\oa=A/k\cdot 1$.  Put
$$
	\deg D\;=({\em {degree\; of\; the\; linear\; map\; }}D)+d  $$
$$
	|D| \;=\deg D-1; \;\;\;|a|=\deg a - 1
$$
Given a tensor $a_1\otimes\cdots\otimes a_N$ in $A^{\otimes N}$, we
will denote it by $(a_1, \ldots, a_n)$. We will write $\eta_j=\sum_{i=1}^{j}|a_i|$
(as in \c{G}).  Put for cochains $D$ and $E$ from $C^*(A,A)$
$$
	(D\smile E)(a_1,\dots,a_{d+e})=(-1)^{\deg E\cdot\eta_d}
	D(a_1,\dots,a_d)\times  
$$
$$
	\times E(a_{d+1},\dots,a_{d+e});
$$
$$
	(D\circ E)(a_1,\dots,a_{d+e-1})=\sum_{j \geq 0}
	(-1)^{|E|\eta_j}
	D(a_1,\dots,a_j,  
$$
$$
	E(a_{j+1},\dots,a_{j+e}),\dots); 
$$
$$
        [D, \; E]= D\circ E - (-1)^{|D||E|}E\circ D
$$
These operations define the graded associative algebra
$(C^*\asa,\deg,\smile)$ and the graded Lie algebra
($C^*\asa$, $|\;\cdot\;|$, $[\;,\;]$) (cf. \c{CE}; \c{Ge}).
Let
$$
	m(a_1,a_2)=(-1)^{\deg a_1}\;a_1 a_2;
$$
this is a 2-cochain of $A$ (not in $C^2$).  Put
$$
	\delta D=[m,D];  
$$
$$
	(\delta D)(a_1,\dots,a_{d+1})=(-1)^{|a_1||D|+|a_1|+1}\times  
$$
$$
	\times a_1 D(a_2,\dots,a_{d+1})+  
$$
$$
	+\sum
_{j=1}^{d}(-1)^{|D|+\eta_j}
		D(a_1,\dots,a_j\;a_{j+1},\dots,a_{d+1})  
$$
$$
	+(-1)^{|D|+\eta_d+1}D(a_1,\dots,a_d)a_{d+1}
$$
The last formula defines a cochain differential on $C^*(A,M)$ for any $M.$
For an element $x$ of $A$, let $\ux$ be the corresponding zero-cochain in $C^*(A,A).$.
By definition
$$
	(\delta\ux)(a)=(-1)^{\deg x}[x,a];
$$
a one-cochain $D$ is a cocycle \iff\ it is a derivation.  One has
$$
	\delta^2=0;\quad\delta(D\smile E)=\delta D\smile E+(-1)^{\deg D}
		D\smile\delta E  
$$
$$
	\delta[D,E]=[\delta D,E]+(-1)^{|D|}\;[D,\delta E]
$$
($\delta^2=0$ follows from $[m,m]=0$).

Thus $C^*\asa$ becomes a complex; the cohomology of this complex 
is $H^*\asa$ or the \hoc cohomology.  The $\smile$ product induces the
Yoneda product on $H^*\asa=\ext_{A\otimes A^0}^*\asa$.  The operation
$[\;,\;]$ is the Gerstenhaber bracket \c{Ge}. 

If $(A, \;\; \partial)$ is a differential graded algebra then one can define the differential $\partial$ acting on $A$ by 
$$
\partial D \;\; = \; [\partial , D]
$$
\begin{definition} Define the differential graded algebra $\E $ as the complex $C^*(A,A)$ with the differential $\delta$ (or $\delta + \partial$ if $(A, \partial )$ is a differential graded algebra), the grading ${deg}$ and the product $\smile .$
\end{definition}

 In Sections 2--4 the only cochains we will be considering will be those from $C^*(A,A).$For \hoc cochains $D$ and $D_i$ define a new \hoc cochain by the following formula of Gerstenhaber (\c{Ge}) and Getzler (\c{G}):

$$
D_0\{D_1, \ldots , D_m\}(a_1, \ldots, a_n) = 
$$
$$
=\sum (-1)^{ \sum_{p=1}^{m}\eta_{i_p}| D_p|}  D_0(a_1, \ldots ,a_{i_1} , D_1 (a_{i_1 + 1}, \ldots ),\ldots ,D_m (a_{i_m + 1}, \ldots ) , \ldots)
$$
\begin{proposition}
 One has
$$
(D\{E_1, \ldots , E_k \})\{F_1, \ldots, F_l \}=\sum (-1)^{\sum _{q \leq i_p}|E_p||F_q|} \times 
$$
$$
\times D\{F_1, \ldots , E_1 \{F_{i_1 +1}, \ldots , \} , \ldots ,  E_k \{F_{i_k +1}, \ldots , \}, \ldots, \}
$$
\end{proposition}
{\ub{Proof}}  Direct computation   $\pf$
 
For a cochain $D$ let $D^{(k)}$ be the following $k$-cochain of $\E$:
$$
D^{(k)}(D_1, \ldots, D_k) = D\{D_1, \ldots, D_k\}
$$
\begin{proposition} \label{etoee}
 The map 
$$
D \mapsto \sum_{k \geq 0} D^{(k)}
$$
is a morphism of differential graded algebras $\E \to \cal{E}^*_{\E}$.
\end{proposition}

\section{The Hochschild chain complex} \label{hochain}

Let $M$ be a bimodule over $A$.
Recall that the \hoc homological complex $(C_*(A,M),b)$ is the following:
$$
	C_n(A,M)=M\otimes\oa^{\otimes n};  
$$
We will always write, as in Section 1, 
$$
(a_0,\dots,a_n) = a_0 \otimes \dots \otimes a_n
$$
For $a_0$ in $M$ and $a_i$ in $A$, $i > 0,$ define
$$
	C_n\asa=A\otimes\oa^{\otimes n};  
$$
$$
	b(a_0,\dots,a_n)=\sum_{j=0}^{n-1}
	(-1)^{\eta_{j+1}+1}(a_0,\dots,a_j\;a_{j+1},\dots,a_n)  
$$
$$
	+(-1)^{(|a_n|+1)(\eta_n+1)+1}
	(a_n a_0,a_1,\dots,a_{n-1}).
$$
We shall denote $C_*(A,A)$ simply by $C_*(A).$ We introduce a grading on $C_*(A,A)$ by the formula
$$\deg(a_0 , \ldots , a_n) = \sum  \deg a_i + n$$
(A rule for remembering the signs:  let $|\epsilon|=1$, $\epsilon^2=0$;
map $C_n(A)$ to $A[\epsilon]\big/[A[\epsilon],A[\epsilon]]$, $a_0\otimes
\cdots\otimes a_n\mapsto a_0\epsilon a_1\epsilon\cdots a_n\epsilon$;
then $b$ becomes $\frac{\partial}{\partial\epsilon}$).

The homology of $C_*(A, M)$ is the \hoc homology $HH_*(A, M)=\tor_*^{A\otimes
A^{\circ}}(A,M)$. We will denote   $HH_*(A, A)$ simply by $HH_*(A).$ If $A$ is a differential graded algebra and
$\partial$ the differential in $A$, one extends $\partial$ to
$C_*(A)$:
$$
	\partial(a_0,\dots,a_n)=\sum_{j=0}^{n}
	(-1)^{\eta_j}(a_0,\dots,\partial a_j,\dots,a_n)
$$
For $a$ in $C_*(A,A)$ and $x$ in $C_*(\E , \E)$ define
\begin{equation} \label{eq:bullet}         
     a\bullet x =  a\bullet _1 x +  a\bullet _2 x
\end{equation}
where
$$
     (a_0 ,\ldots , a_n) \bullet_1  (D_0 , \ldots ,      D_m) = \sum (-1)^{(\eta _{n+1} - \eta_1) \deg D_0 + \sum_{p=1}^{m}{(\eta _{n+1} - \eta_{i_p + 1})| D_p|}}\times 
$$
\begin{equation} \label{eq:bullet1}
 \times (a_0 D_0(a_1, a_2, \ldots ) ,\ldots ,a_{i_1} , D_1 (a_{i_1 + 1}, \ldots ),\ldots ,D_m (a_{i_m + 1}, \ldots ) , \ldots)
\end{equation}
\begin{equation} \label{eq:bullet2}
     (a_0 ,\ldots ,a_n) \bullet_2  (D_0 , \ldots ,      D_m) = 
\end{equation}
$$
=\sum_{q \leq n+1} (-1)^{ (\eta_{n+1}-\eta_q) \eta _q + (\eta _q - \eta _{i_0 + 1}) \deg D_0 + \sum_{p=1}^{m-1}{(\eta _q - \eta_{i_p + 1})| D_p|} + \eta _{n+1} |D_m|}
$$
$$
 \times (-1)^{ |D_m|(\sum_{p \geq 0}|D_p|+1)} \times 
$$
$$
\times (D_m(a_q, \ldots , a_n, a_0, \ldots, a_{i_0}) D_0(a_{i_0 + 1}, \ldots ) ,\ldots , a_{i_1} , 
$$
$$
, D_1 (a_{i_1 + 1}, \ldots ), \ldots , D_{m-1} (a_{i_{m-1} + 1}, \ldots ) , \ldots )
$$
The sum in (\ref{eq:bullet1}) is taken over all $q, \ i_0, \ldots , \ i_{m-1}$ for which $a_0$ is inside $D_m.$
\begin{th}
 The map
$$
      \bullet : C_*(A,A) \otimes C_*(\E , \E) \to C_*(A,A)
$$
is a morphism of complexes, i.e.
$$
       b(a \bullet x) = (ba) \bullet x + (-1)^{\deg a} a \bullet (b + \delta)x
$$
\end{th}
{\ub{Proof}} Direct computation. $\pf$

Consider some examples of the product  $\bullet$. When all $D_i$ are zero -cochains then one gets the shuffle product from \c{CE}. On the other hand, 
put, for a $d$-cochain $D$ and for $a$ in $C_*(A,A)$, $i_D(a) = (-1)^{\deg a\;\deg D}(a \bullet D);$
$$
	i_D(a_0,\ldots,a_n)=(-1)^{\deg a_0 \cdot \deg D_0}
	(a_0 D(a_1,\dots,a_d),\quad a_{d+1},\dots,a_n)
$$

One gets the cap product from \c{CE}:
$$
	C^d(A,A)\otimes C_n(A,A)\to C_{n-d}(A,A)
$$
Now, for a $d$-cochain $D$ and for $a$ in $C_*(A,A)$, let $L_D (a) = (-1)^{|D| \cdot \deg a} a \bullet (1, D);$ one has
$$
L_D(a_0, \ldots,a_n)= 
\sum_{q=0}^{n-d+2} (-1)^{(\eta _{n+1} - \eta _q) \eta _q + |D|} 
(D(a_q , \ldots, a_0 , \ldots), \ldots, a_{q-1}) +
$$
$$ 
+\sum_{k=1}^{n-d}(-1)^{(\eta _{k+1}+1) |D|}(a_0, \ldots, a_k , D(a_{k+1}, \ldots ), \ldots, )
$$
One has
$$
     [b, L_D] = -L_{\delta D}
$$
Now construct the product

\begin{equation}  \label{eq:bullet0}
      \bullet : C_*(\E, \E) \otimes C_*(\E , \E) \to C_*(\E , \E)
\end{equation}
 by composing the one for the algebra $\E$ with the map $\E \to \cal{E}^*_{\E}$ from Proposition \ref{etoee}.

The product (\ref{eq:bullet0}) is given by formulas (\ref{eq:bullet1}), (\ref{eq:bullet2}) where $a_i$ are now viewed as \hoc cochains and $D(a_1, \ldots, a_d)$ is replaced by $D\{a_1, \ldots, a_d\}$.
\begin{proposition}
The product (\ref{eq:bullet0}) is homotopically associative.
\end{proposition}
We omit the proof.

\section{The periodic cyclic complex} \label{ccper}
In the situation of Sect. \ref{hochain} define
$$
	B(a_0,\dots,a_n)=\sum_{j=0}^{n}
	(-1)^{(\eta_{n+1}-\eta_j)\eta_j}(1,a_j,\dots,a_n, 
	a_0,\dots,a_{j-1})
$$
Then $B^2=b^2=Bb+bB=0$; denote
$$
	CC_n^{\per}(A)=\prod_{i \equiv n (mod 2)}
	A\otimes\oa^{\otimes n}
$$
with the differential $b+B+\partial$ (for the differential graded
algebra $(A,\partial)$).  Put 
\begin{equation}  \label{eq:bullet3}
	(a_0, \ldots, a_n) \bullet _3 (D_1, \ldots, D_m) = 
\end{equation}
$$
	=\sum _{0\leq p\leq m}
	(-1)^{\sum_{r<p} (|D_r|+1)\sum_{r\geq p}	       (|D_r|+1)+\sum|D_r|}\times 
$$
$$
	\times\sum _ 	{0\leq j\leq n}			
	(-1)^{(\eta_{n+1}-\eta_j)\eta_j}\lambda'
	(D_p,\dots,D_m,D_0,\dots,D_{p-1})  
	(1,a_j,\dots,a_n,\;a_0,\dots,a_{j-1})
$$
Here

$$
	\lambda(D_1,\dots,D_m)(a_0,\dots,a_n)= 
 $$
$$
	=\sum_{j_1>0}
	(-1)^{|D_1|(\eta_{n+1}-\eta_{j_1})+\dots+|D_m|(\eta_{n+1}-\eta_{j_m})}\times 
$$
$$
	\times(a_0,\dots,D_1(a_{j_1}\dots),\dots,D_m(a_{j_m},\dots),\dots);  
$$

 $\lambda'$ is the sum of all those terms in $\lambda$ 
which contain $D_0(\dots)$ after $a_0$.

\begin{th}  \label {bulletnaccper}
 The map
$$
      \bullet : \ccp(A) \otimes \ccp(\E) \to \ccp(A)
$$
is a morphism of complexes.
\end{th}
{\ub{Proof}} Direct computation.  $\pf$

For a $d$-cochain $D$ and for $a$ in $\ccp(A)$, let $I_D(a) = (-1)^{\deg D \cdot \deg a}(a \bullet D);$ $L_D (a) = (-1)^{|D| \cdot \deg a} a \bullet (1, D);$ one has

$$
	[B+b,\;I_D]-I_{\delta D}=L_D.  \pf 
$$

This is the homotopy formula of Rinehart \c{R}.

One can define the product
$$
      \bullet : \ccp(\E) \otimes \ccp(\E) \to \ccp(\E)
$$
exactly in the same way as in the end of Section 2.

\begin{remark}  \label{joneshood}
 Considering $CC^{per}_*(A)=CC^{per}_*({\cal{E}}_A^0)$ as a part of $CC^{per}_*(\E)$ and restricting the product $\bullet$ to it, we get the Hood-Jones product on $CC^{per}_*(A).$ When $A$ is commutative then $CC^{per}_*(A)$ is a subcomplex of $CC^{per}_*(\E)$ and one gets a product on the complex $CC^{per}_*(A).$
\end{remark}

\section{The Fukaya category}  \label{fukaya}
\begin{equation} \label{eq:ca}
{\cal{C}}^n(A) = \bigoplus_{i+j=n} C_{-i}({\cal{E}}^*,{\cal{E}}^*)^j
\end{equation}
\begin{equation}  \label{eq:cainf}
{\cal{C}}^n_{\infty}(A) = \prod_{i+j=n} C_{-i}({\cal{E}}^*,{\cal{E}}^*)^j
\end{equation}
In this Section we extend the results of Section \ref{hochain} by ``moving away from the diagonal''. For any automorphism $\alpha$ of an algebra $A$ we define the twisted Hochschild complex $C_*(A,A_{\alpha})$ (related to noncommutative geometry of fixed points of $\alpha$) and construct in these terms a homotopically associative category for which ${\cal{C}}^*(A)$ or ${\cal{C}}^*_{\infty}(A)$ is the ring of endomorphisms of the identity.

Let $\alpha,$ $\beta$ be two automorphisms of $A.$ Define the new bimodule $_{\alpha} A _{\beta}$ over $A$ as follows: $ _{\alpha} A _{\beta} = A$ as $k$-modules and $a\cdot m \cdot b = \alpha(a) m \beta (b)$ for $a,m,b$ in $A.$
Consider the chain complex $C_*(A,\;_{\alpha} A _{\beta})$ as in Section 2 and the cochain complex $_{\alpha}{\cal{E}}^*_{\beta} = C^*(A,\;_{\alpha} A _{\beta})$ as in Section 1. We shall write $A_{\alpha} = _{{id}} A _{\alpha}$. The cup product on \hoc cochains is well defined as a morphism $_{\alpha}\cal{E}^*_{\beta} \otimes _{\beta}\cal{E}^*_{\gamma} \rightarrow_{\alpha}\cal{E}^*_{\gamma}$. 

Given an automorphism $\alpha$ of $A$, one can define its action on a \hoc 
cochain from $\E$ in two ways:

$$
(\alpha D)(a_1, \ldots, a_n) = \alpha (D(a_1, \ldots, a_n));
$$
$$ 
(D \alpha)(a_1, \ldots, a_n) =  D(\alpha a_1, \ldots,\alpha a_n)
$$
Both yield morphisms of complexes $\E \rightarrow _{\alpha}\cal{E}_{\alpha}$.
We define an $\E$-bimodule structure on $_{\alpha}\cal{E}^*_{\beta} $ as follows:
$$
D\cdot M \cdot E = D\alpha \smile M \smile \beta E 
$$
for $D,E \in \E$ and $M \in $$_{\alpha}\cal{E}^*_{\beta} $.

\begin{th}
The pairings (\ref{eq: oper}, \ref{eq: oper1}) can be extended to the homotopically associative natural morphisms of complexes
$$ \bullet :C_*(A,A_{\alpha}) \otimes C_*(\E, _{\alpha}{\cal{E}}^*_{\beta}) \rightarrow  C_*(A, A_{\beta})
$$
$$
\bullet :C_*(\E,\; _{\alpha}{\cal{E}}^*_{\beta}) \otimes C_*(\E, \;_{\beta}{\cal{E}}^*_{\gamma}) \rightarrow C_*(\E, \;_{\alpha}\cal{E}^*_{\gamma})
$$
\end{th}
{\ub {Proof}}
As in Section 2, put for $a \in C_*(A,A_{\alpha})$ and $x \in C_*(\E, \; _{\alpha}\cal{E}_{\beta})$
$$         
     a\bullet x =  a\bullet _1 x +  a\bullet _2 x
$$
where
$$
     (a_0 ,\ldots , a_n) \bullet_1  (D_0 , \ldots ,      D_m) = 
$$
$$
\sum \pm 
 (a_0 \cdot D_0(a_1, \ldots ) ,\ldots ,a_{i_1} , D_1 (a_{i_1 + 1}, \ldots ),\ldots ,D_m (a_{i_m + 1}, \ldots ) , \ldots)
$$
$$
     (a_0 ,\ldots ,a_n) \bullet_2  (D_0 , \ldots ,      D_m) = 
$$
$$
=\sum_{q \leq n+1} \pm (D_m (a_q, \ldots , a_n, a_0, \alpha a_1, \ldots, \alpha a_{i_0})\cdot D_0(a_{i_0 + 1}, \ldots ) ,\ldots , a_{i_1} , 
$$
$$
, D_1 (a_{i_1 + 1}, \ldots ), \ldots , D_{m-1} (a_{i_{m-1} + 1}, \ldots ) , \ldots )
$$

For $A \in  C_*(\E, \; _{\alpha}\cal{E}_{\beta})$ and $x \in  C_*(\E, \; _{\beta}\cal{E}_{\gamma})$

$$         
     A\bullet x =  A\bullet _1 x +  A\bullet _2 x
$$
where
$$
     (A_0 ,\ldots , A_n) \bullet_1  (D_0 , \ldots ,      D_m) =
$$
$$
\sum \pm 
 (A_0 \smile D_0\{A_1, A_2, \ldots \} ,\ldots ,A_{i_1} , D_1 \{A_{i_1 + 1}, $$
$$
\ldots \},\ldots ,D_m \{A_{i_m + 1}, \ldots \} , \ldots)
$$
$$
     (A_0 ,\ldots ,A_n) \bullet_2  (D_0 , \ldots ,      D_m) =
$$
$$
=\sum_{q \leq n+1} \pm (D_m\{\{A_q, \ldots , A_n, A_0, A_1, \ldots, A_{i_0}\}\} \smile D_0\{A_{i_0 + 1}, \ldots \} ,\ldots , A_{i_1} ,
$$
$$
, D_1 \{A_{i_1 + 1}, \ldots \}, \ldots , D_{m-1} \{A_{i_{m-1} + 1}, \ldots \} , \ldots )
$$
and 
$$
D_m\{\{A_q, \ldots , A_n, A_0, A_1, \ldots, A_{i_0}\}\}(a_1, \ldots, a_n) =
$$
$$
= \sum \pm D_m (\alpha a_1, \dots, A_q (\alpha a_j, \dots, ), \dots,
$$
$$
, \alpha a_t,  
 A_0 (a_{t+1}, \dots, ),\beta a_p, \dots, \beta A_1(a_r, \dots ), \beta a_s, \dots )
$$

One checks that these maps are homotopically associative morphisms of complexes.
$\pf$

One gets the cohomology groups $\cal{H}^*(\alpha, \beta)$ and  $\cal{H}^*_{\infty}(\alpha, \beta)$ which form a category ${\cal H}^*_A$ or  ${\cal H}^*_{A, \infty}.$

\section{Examples} \label{examples}
\begin{th}   \label {diffoponforms}
 Let $A = k[X]$ be the algebra of regular functions on an affine nonsingular algebraic variety $X$ over a field $k$ of characteristic zero. Then
$$
{\cal{H}}^*(A) \rightisoarrow {\cal{D}}(\Omega^*_{X/k})^{op}
$$
(the ring opposite to the ring of all differential operators on $\Omega$);
$$
{\cal{H}}^*_{\infty}(A) \rightisoarrow{End}(\Omega^*_{X/k})^{op}
$$
\end{th}
Note that $HH_*(A) \simeq \Omega^*_{X/k}$; the $\bullet$ action of ${\cal{H}}^*$ on $HH^*$ is the obvious one.

{\ub{Proof of \ref{diffoponforms}}} 
\begin{lemma}   \label{quisofdga}
 Let ${\cal{E}}^*,$ ${\cal{F}}^*$ be two differential graded algebras together with a quasi-isomorphism $f: {\cal{E}}^* \to {\cal{F}}^*.$ Then $f$ induces an isomorphism
$$
\bigoplus_{i+j=n} C_{-i}({\cal{E}}^*,{\cal{E}}^*)^j \rightisoarrow \bigoplus_{i+j=n} C_{-i}({\cal{F}}^*,{\cal{F}}^*)^j
$$
\end{lemma}
{\ub{Proof}} The corresponding spectral sequence converges (this is not true if one replaces $\bigoplus$ by $\prod$).     $\pf$

Note that by \c{HKR} $\E$ is quasi-isomorphic to ${\cal{F}}^* = \Gamma (X, \bigwedge ^*TX)$ as differential graded algebras. Locally, ${\cal{F}}^*_U = {\cal{O}}_U \otimes \bigwedge ^*(\partial_ {x_1}, \ldots, \partial_ {x_n});$ 
$HH_*({\cal{F}}^*_U) = \Omega^*_{U/k} \otimes \bigwedge ^*(\partial_{x_1}, \ldots, \partial_ {x_n}) \otimes S^*(d\partial_{x_1}, \ldots, d\partial_ {x_n})$ (classes $d\partial_{x_i}$ are represented by cycles $(1, d\partial_{x_i})$). Now, looking at the $\bullet$ action on $HH^*({\cal O}_U) \simeq \Omega ^*_{U/k}$, one sees that $\Omega ^*_{U/k}$ acts by multiplication, $\partial _{x_i}$ acts by contraction $i_{\partial _{x_i}}$ and $d\partial _{x_i}$ acts by Lie derivative $L_{\partial _{x_i}}.$ Whence the local isomorphism
$$
{\cal{H}}^*({\cal{O}}_U) \rightisoarrow {\cal{D}}(\Omega^*_{U/k})^{op};
$$
this statement may be globalized exactly as in \c{HKR}
Now let us compute ${\cal{H}}^*_{\infty}(A)$. We have to use the other spectral sequence whose $E_1$ term is $HH_*(\E)$ if one forgets about the differential $\delta$. But, as a graded algebra, $\E$ is just $T^*(A') \widehat{\otimes} A$ by which we mean the space of all multilinear maps from $A$ to $A$. One can show that 
$$
   E_2^{-p, n} = HC^{n-1}(A)  \widehat{\otimes} \Omega ^p +  HC^{n-2}(A)  \widehat{\otimes} \Omega ^{p+1}
$$
$$
   E_3^{-p, n} = (\Ker S) ^{n-1}  \widehat{\otimes} \Omega ^p +  (\Coker S)^{n-2}  \widehat{\otimes} \Omega ^{p+1}
$$
where $S : HC^n \to HC^{n+2}$ is the Bott operator (\c{L}).
Recall the Gysin exact sequence (\c{C}, \c{T}, \c{L}):
$$
    0 \longrightarrow (\Coker S)^{n} \longrightarrow HH^n(A) \longrightarrow (\Ker S)^{n-1} \longrightarrow 0 
$$
It is not hard to construct for any element of $HH^*(A) \widehat{\otimes} \Omega ^*_{X/k}$ a corresponding cycle of the complex ${\cal{C}}^*$. This shows that the spectral sequence degenerates at $E_3$ term and that ${\cal{H}}^*_{\infty}(A) \rightisoarrow HH^*(A) \widehat{\otimes} \Omega ^*_{X/k},$     $\pf$
 
Let $A$ be a one-dimensional affine space and $A=k[V] = k[x].$ We will compute the full subcategory of ${\cal H}^*_A$ whose objects are dilations $x \mapsto \alpha x.$ We view ${\cal D}={\cal D}(\Omega^*(V))$ as a graded algebra $k[x,\zeta, \frac{\partial}{\partial x}, \frac{\partial}{\partial \zeta}]$ where $\zeta = dx .$
\begin{proposition} 
$$\H ^*(1,1) \rightisoarrow {\cal D}^{op};$$
$$\H ^*(\alpha, \beta) \rightisoarrow k, \;\;\;\; \alpha \neq 1,\;\;\; \beta \neq 1 \;\; ;$$
$$\H ^*(\alpha, 1) \rightisoarrow k \otimes_{k[x, \partial_{\zeta}]} {\cal D}^{op};$$
$$\H ^*(1, \alpha) \rightisoarrow {\cal D}^{op}\otimes_{k[x, \partial_{\zeta}]} k .$$
\end{proposition}

We denote by $V^{\alpha}$ the set of fixed points of $\alpha$ 
and consider ${\cal D}(\Omega^*(V^{\alpha}))$ as the algebra of functions on the supermanifold $TV^{\alpha}.$ 
Let us identify $TV$ with $T^*V;$ consider the maps 

$$
TV^{\alpha} \stackrel {j_{\alpha}}{\longleftarrow} T_{V^{\alpha} \cap V^{\beta}}V \stackrel {j_{\beta}}{\longrightarrow} TV^{\beta}
$$
where $j_{\alpha}$ is the composition
$$ T_{V^{\alpha} \cap V^{\beta}}V \stackrel {i}{\longrightarrow} T_{V^{\alpha}}V \stackrel {\pi}{\longrightarrow} TV.
$$
Then 
$$
\H^*  (\alpha, \alpha)={\cal D}(\Omega ^*(V^{\alpha}))^{op};
$$
$$
\H^* (\alpha, \beta)=j_{\alpha *}j_{\beta}^*{\cal D}(\Omega ^*(V^{\beta}))
$$
which is a left $\H (\alpha, \alpha)$-module and a right $\H (\alpha, \alpha)$-module.

Denote by $e_{\alpha \beta}$ the generators of $\H^0(\alpha, \beta)$ for $\alpha \neq 1$ and by $D_{\alpha \beta}$ the generators of $\H^1 (\alpha, \beta)$ for $\alpha = 1$ or $\beta = 1.$

\begin{proposition}
One has
$$
e_{\alpha \beta} \bullet e_{ \beta \gamma} = e_{\alpha \gamma};\;\;\;\; 
e_{\alpha \beta} \bullet D_{ \beta \gamma} = D_{\alpha \gamma};\;\;\;\;
$$
$$ 
D_{\alpha \beta} \bullet e_{ \beta \gamma} = D_{\alpha \gamma};\;\;\;\; 
D_{\alpha \beta} \bullet D_{ \beta \gamma} = 0; 
$$
the right, resp.left ${\cal D}^{op}$-module structure on $\H^*(1,\alpha),$
resp.  $\H^*(\alpha,1),$ is the obvious one.
\end{proposition}
\begin{proposition}
One has
$$
H_*(A,A)=\Omega^*(V); \;\;\; H_*(A,A_{\alpha})=k, \;\;\;\; \alpha \neq 1.
$$
Let $1_{\alpha}$ be the generator of $ H_0(A,A_{\alpha}).$ Then for $\omega \in \Omega^*(V)$ $\omega \bullet D_{1\alpha} = (i_{\frac{\partial}{\partial x}}\omega)(0);$ $1_{\alpha} \bullet D_{\alpha \beta} = 1_{\beta}$ if $\alpha, \beta \neq 1;$ $1_{\alpha} \bullet D_{\alpha 1}=0.$
\end{proposition}

In other words, $H_*(A,A_{\alpha})=\Omega^*(V^{\alpha})$ and one has
$$
\omega \bullet D_{\alpha \beta} = j_{\beta *}j_{\alpha}^*\omega .
$$
Here $j_{\alpha}^*$ is the usual restriction; $j_{\beta *}$ is either zero when $j_{\beta }$ is a closed embedding or the Berezin integral when $j_{\beta }$ is a projection along the odd fiber.

{\ub {Proof of the Propositions.} The \hoc homology $_{\alpha}H^*_{\beta } \stackrel{def}{=}H^*(A,_{\alpha}A_{\beta })$ is computed by the complex $k[x] \stackrel{\alpha - \beta}{\rightarrow} k[x].\;\;$So for $\alpha \neq \beta$  $_{\alpha}H^1_{\beta } \simeq k\;\;\;$ (we denote the generator by $D_{\alpha \beta}$)and $ _{\alpha}H^i_{\beta } =0 $ if $i \neq 1.$ Consider the spectral sequence with the $E_2$ term $H_*(H^*, _{\alpha}H^*_{\beta } )$ where $H^*=_{1}H^*_{1 }.$ When $\alpha = \beta,$ we get the ring $H^*=k[x,\partial _x]$ acting on itself; $x\cdot a=xa;\;\;a\cdot x =\alpha a x;\;\;\partial_x \cdot a=\alpha \partial_x a;\;\;a\cdot \partial_x =a \partial_x .$So for $\alpha = \beta \neq 1$ $E_2 = k$ and the spectral sequence collapses. 

When $\alpha \neq \beta$ then the $E_2$ term is $k[dx, d\partial_x] \cdot D_{\alpha \beta}$ where $dx, d\partial _x$ are the generators from $H_1(H^*,k).$ We claim that the spectral sequence degenerates at $E_2$ if $\alpha =1$ or $\beta = 1$ and that $E_3 = k$ if $\alpha,\;\; \beta \neq 1.$ Indeed, if $\alpha,\;\; \beta \neq 1$ then 
$$
d_2(D_{\alpha \beta} \cdot dx \cdot (d\partial _x)^n = (\alpha - \beta)(\alpha - 1)(\beta - 1) (D_{\alpha \beta} \cdot (d\partial _x)^n;
$$
so all terms vanish in $E_3$ except $D_{\alpha \beta} \cdot dx.$ This terms always survives in $E_{\infty}$ if $\alpha \neq \beta.$ It is represented by the cycle 
$$
(D_{\alpha \beta},x) + \frac{\beta - 1}{\alpha - \beta} \cdot 1.
$$
In case when $\alpha =1,\;\;$ resp. $\beta = 1,\;\;$ the element above is equal to $x \bullet D_{\alpha \beta},\;\;$ resp.  $ D_{\alpha \beta}  \bullet x.\;\;$
Therefore $D_{\alpha \beta} \neq 0;  \;$ it generates the ${\cal{D}}^{op}$-module ${\cal {H}}^*(\alpha, \beta)$ and is annihilated by $x$ and $\frac{\partial}{\partial \zeta}.$ This determines this module uniquely.
The rest of the statements of the Propositions are obtained by direct computations (note that $H_*(A, A_{\alpha})$ is the homology of the complex $k[x] \stackrel {\alpha}{\rightarrow} k[x]$). $\pf$

Now let us turn to the case when $A$ is a deformed algebra on a symplectic manifold. Let $*$ be a star product on a symplectic manifold $M.$ Denote by $\ah (M)$ the algebra $(C^{\infty }(M), *).$ 
\begin{th}
One has 
$${\H}^*(\ah (M)[\hbar ^{-1}]) \rightisoarrow H^*(M^{S^1}, \C [[\hbar, \hbar ^{-1}])$$
\end{th}
{\ub {Proof.}} 
\begin{lemma}
The differential graded algebra ${\cal {E}}^* _{\ah (M)[\hbar ^{-1}]}$ is equivalent to  the algebra $(\Omega ^*(M)[[\hbar, \hbar ^{-1}], d).$
\end{lemma}
By them being equivalent we mean that they are connected by a chain of quasi-isomorphisms of differential graded algebras. The Theorem follows from the Lemma because of Lemma \ref{quisofdga}. To prove the Lemma, consider the Weyl bundle $\W$ with a Fedosov connection $\Nabla ^F$ (\c{Fe}). Consider the double complex $\Omega^* (M, C^*(\W))$ with the differential $\Nabla ^F.$ This is a differential graded algebra which, after inverting $\hbar,$ is quasi-isomorphic to both algebras in the statement of the Lemma. Indeed, on the one hand, for every $p$ the complex  
$\Omega^* (M, C^p({\W}))[\hbar ^{-1}]$ has only zero-dimensional cohomology group which is $C^p (\ah (M)[\hbar ^{-1}]);$ on the other hand, the embedding $\C [[\hbar, \hbar ^{-1}] \rightarrow C^*(\ah (M)[\hbar ^{-1}])$ is a quasi-isomorphism invariant under all derivations of $\ah(M)$ (which can be seen by filtering the complex by powers of $\hbar$ and considering the corresponding spectral sequence).                $\pf$

Next we will compute ``the Fukaya subcategory'' of a deformed algebra of functions on $M = S^2.$ We will consider the deformation of the algebra of algebraic functions. Let $A = U({\frak{sl}} _2) / (\Delta - c^2)$ where $\Delta$ is the Casimir operator. Consider the subgroup $S^1$ of $SU_2$ acting on $A;$ 
$$
X _{\pm} \mapsto \alpha ^{\pm 1} X_{\pm} ; H \mapsto H
$$
where $\alpha \in S^1 .$

Let $_{\alpha} \Omega _{\beta} (M) = \{ l: S^1 \rightarrow M : \alpha l(1) = \beta l(1)\}.$
\begin{proposition}
One has 
$$
{\cal H}^*_A(\alpha, \beta)= H^{*-d_{\alpha \beta}}(_{\alpha} \Omega _{\beta} (M))
$$
where $d_{\alpha \beta}={\em codim} \{x | \alpha x = \beta x \};$
\begin{equation} \label{eq: haha}
H_*(A,A_{\alpha}) = H^{2-*}(M), \;\;\;\alpha = 1;
\end{equation}
\begin{equation} \label{eq: haha1}
H_*(A,A_{\alpha}) = \bigoplus_{\alpha x = x}{\C}, \;\;\; \alpha \neq 1.
\end{equation}
\end{proposition}
{\ub{Proof.}} Recall that if $V$ is a $U({\g})$-bimodule then $H^*(U({\g}), V)$ is computed by the Koszul complex $(C_*({\g}, V), \partial)$ where the $\g$-module structure on $V$ is given by $X \cdot v = Xv-vX.$ Consider the operator $D: C_* \rightarrow C_{*+1},$
$$
D(v \otimes (X_1 \wedge \ldots \wedge X_n)) = \sum _{i} (e_i v + v e_i) \otimes (e^i \wedge X_1 \ldots \wedge X_n)
$$
where $(e_i, e^i)$ are dual bases with respect to the Killing form. The mixed complex 
$$
(R_n = \bigoplus C_{n-2i}({\g}, A \otimes A); \;\; \partial + D)
$$
is a free bimodule resolution of $A$ (the subcomplex of the standard bar resolution). Computing the spectral sequences for $R_*\otimes _{A \otimes A^o} A_{\alpha}$ or ${\em Hom}_{A \otimes A^o}(R_*, _{\alpha} A _{\beta}),$ we see that they degenerate at $E_2.$ One sees that the homology $H_*$ is given by (\ref{eq: haha}, \ref{eq: haha1}) and the cohomology is given by 
$$
H^*(A, _{\alpha} A _{\beta}) = H^*(S^2),\;\;\; \alpha = \beta ;
$$
$$
H^*(A, _{\alpha} A _{\beta}) = \bigoplus _{\{x \;|\; \alpha x= \beta x\}} {\C}[-2], \;\;\; \alpha \neq \beta.
$$
The spectral sequence with the $E_2$ term $H_*(H^*, _{\alpha} H^*_{\beta})$ degenerates at $E_2$ (there is no place for a nonzero differential).
$\;\;\pf$

\end{document}